\documentclass[12pt,a4paper,psamsfonts]{amsart}
\usepackage{amssymb,amscd,amsxtra,calc}
\usepackage{cmmib57}

\theoremstyle{plain}
    \newtheorem{thm}{Theorem}[section] %
    \renewcommand{\thethm}%
    {\arabic{section}.\arabic{thm}}
    \newtheorem{claim}[thm]{Claim}
    
    \newtheorem{corollary}[thm]{Corollary}
    \newtheorem{example}[thm]{Example}

    \newtheorem{question}[thm]{Question}
    \newtheorem{theorem}[thm]{Theorem}

\theoremstyle{definition}
    \newtheorem{definition}[thm]{Definition}

    \newtheorem{remark}[thm]{Remark}
\theoremstyle{remark}

    \newtheorem{setup}[thm]{}

\newcommand{\BCC}{\mathbb{C}}

\newcommand{\BPP}{\mathbb{P}}
\newcommand{\BQQ}{\mathbb{Q}}
\newcommand{\BRR}{\mathbb{R}}
\newcommand{\BZZ}{\mathbb{Z}}

\newcommand{\SO}{\mathcal{O}}

\newcommand{\alb}{\operatorname{alb}}
\newcommand{\Aut}{\operatorname{Aut}}

\newcommand{\Exc}{\operatorname{Exc}}
\newcommand{\Gal}{\operatorname{Gal}}
\newcommand{\group}{\operatorname{group}}
\newcommand{\GL}{\operatorname{GL}}

\newcommand{\id}{\operatorname{id}}
\newcommand{\Imm}{\operatorname{Im}}
\newcommand{\Ker}{\operatorname{Ker}}

\newcommand{\NS}{\operatorname{NS}}

\newcommand{\pt}{\operatorname{pt}}
\newcommand{\SL}{\operatorname{SL}}
\newcommand{\torsion}{\operatorname{torsion}}
\newcommand{\variety}{\operatorname{variety}}

\newcommand{\Alb}{\operatorname{Alb}} 
\newcommand{\Spec}{\operatorname{\mathit{Spec}}}

\newcommand{\ratmap}
{{\,\cdot\negmedspace\cdot\negmedspace\cdot\negmedspace\to\,}}

\newcommand{\Lehmer}{\mathrm{Lehmer}}
\newcommand{\rank}{\mathrm{rank}}

\begin{document}

\title[Automorphism groups of projective varieties]
{Automorphism groups of positive entropy on minimal projective varieties}
\author{De-Qi Zhang}
\address
{
\textsc{Department of Mathematics} \endgraf
\textsc{National University of Singapore, 2 Science Drive 2,
Singapore 117543}}
\email{matzdq@nus.edu.sg}

\begin{abstract}
We determine the geometric structure of a minimal projective threefold
having two `independent and commutative' automorphisms of positive topological entropy, and
generalize this result to higher-dimensional smooth minimal pairs $(X, G)$.
As a consequence, we give an effective lower bound for the first dynamical degree
of these automorphisms of $X$ fitting the `boundary case'.
\end{abstract}

\subjclass[2000]{32H50, 14J50, 32M05, 14J32}
\keywords{automorphism group, topological entropy, Calabi-Yau variety}

\thanks{The author is supported by an ARF of NUS}

\maketitle


\section{Introduction}

A normal projective variety $X$ with only
$\BQQ$-factorial terminal singularities is called {\it minimal} if the canonical divisor
$K_X$ is nef.
For $G \le \Aut(X)$,
the representation $G | \NS_{\BCC}(X)$ on the
complexified N\'eron-Severi group $\NS_{\BCC}(X) := \NS(X) \otimes_{\BZZ} \BCC$,
is {\it $Z$-connected}
if its Zariski-closure in $\GL(\NS_{\BCC}(X))$ is connected in the Zariski topology.
We denote by $q(X) := h^1(X, \SO_X)$ the {\it irregularity} of $X$.
A birational morphism $\sigma : X \to X'$ is {\it crepant} if $K_X = \sigma^*K_{X'}$.
We refer to \cite[Definition 2.34]{KM} for the definition of terminal or canonical singularity,
and to \cite{DS} for the definitions of dynamical degrees and (topological) entropy.

In this paper, we prove Theorem \ref{ThA} below and its generalization in Theorem \ref{ThB}.

\begin{theorem}\label{ThA}
Let $G \le \Aut(X)$ be an automorphism group on a minimal
projective threefold $X$
with the representation $G | \NS_{\BCC}(X)$ solvable
and $Z$-connected. Then we have:
\begin{itemize}
\item[(1)]
The null subset $N(G) := \{g \in G \, | \, g$ is of null entropy$\}$ is
a normal subgroup of $G$ such that $G/N(G) \cong \BZZ^{\oplus r}$ for some
$r = r(G) \le \dim X - 1 = 2$.
\item[(2)]
Suppose that $r = 2$. Then there are a $G$-equivariant $($birational$)$ crepant morphism $X \to X'$
and a $G$-equivariant finite Galois cover $\tau: A \to X'$ \'etale in codimension $1$
for a $3$-dimensional abelian variety $A$, with $\deg \tau \ge 2$ occurring only when
$q(X) = 0$.
\item[(3)]
Suppose that $r = 2$ and the identity component $\Aut(X)_0$ of $\Aut(X)$ is trivial.
Then $|N(G)| < \infty$.
\end{itemize}
\end{theorem}

Theorem \ref{ThA} (1) above holds in any dimension $n \ge 2$ and even for any K\"ahler manifold $X$
(with $\NS_{\BCC}(X)$ replaced by $H^{1,1}(X)$), i.e.,
$r(G) \le n -1$ (cf.~\cite[Theorem 1.2, Remark 1.3]{Z4}; \cite{Ki}).
This paper is to prove Theorem \ref{ThA} (2) (and (3)) about the geometric structure
of $X$ when $r(G) = n - 1 = 2$, and generalize it to higher dimensions (cf.~Theorem \ref{ThB}).

As a consequence, we give an effective lower bound for the first dynamical degree
$d_1(g)$ of $g \in G$ when $r(G) = n-1$ (cf.~Corollary \ref{cThB}).

Inspired by the main result of Dinh-Sibony \cite[Theorem 1]{DS}, it was asked
in \cite[Question 2.18]{Z4} whether $r(G) = n-1$ and (the identity component) $\Aut_0(X) = (1)$
imply that the null subset $N(G) \subset G$ (as in Theorem \ref{ThA}) is a finite group.
Theorem \ref{ThA} (3)
answers this question in the affirmative for minimal threefolds.
When $G$ is abelian with $r(G) = n-1$, the finiteness of $N(G)$ is confirmed for all K\"ahler manifolds
and in any dimension by Dinh and Sibony in \cite[Theorem 1]{DS}.

The motivation to consider the assumption in \ref{ThA} is due to
the Tits type result: for any $G \le \Aut(X)$ on a K\"ahler (resp. projective) manifold $X$,
either $G \ge \BZZ * \BZZ$ (non-abelian free group of rank two), or $G$ has a finite-index subgroup $G_1$
such that $G_1|H^{1,1}(X)$ (resp. $G_1 | \NS_{\BCC}(X)$) is solvable and $Z$-connected
(see also \cite[Theorem 1.1, Remark 1.3]{Z4}).

Theorem \ref{ThA} has a special consequence below, which
explains why it seemed hard to construct
two `independent and commutative' automorphisms of positive topological entropy on a Calabi-Yau threefold
except perhaps by descending those on an abelian variety to its quotients.
An invertible linear map is of {\it positive entropy}
if its spectral radius ($=$ the maximum of modulus of eigenvalues)
is greater than one.

\begin{corollary}\label{CorB}
Let $G \le \Aut(X)$ be an automorphism group on a minimal
projective threefold $X$ such that the representation $G^* := G | \NS_{\BCC}(X)
\cong \BZZ^{\oplus 2}$ and
every non-trivial element of $G^*$ is of positive entropy.
Then there are a $G$-equivariant $($birational$)$ crepant morphism $X \to X'$
and a $G$-equivariant finite Galois covering $A \to X'$ \'etale in codimension $1$
for a $3$-dimensional abelian variety $A$.
\end{corollary}

\begin{remark}\label{rThA}
\begin{itemize}
\item[(1)]
The triviality of the identity component $\Aut_0(X)$ of $\Aut(X)$ is a necessary hypothesis
for the finiteness of $N(G)$ in Theorem \ref{ThA} (3). Indeed, a bigger group
$G_1 := G . \Aut_0(X)$ satisfies also the assumption of Theorem \ref{ThA} with
$r(G_1) = r(G)$ and $N(G_1) \ge \Aut_0(X)$.

\item[(2)]
By the remark above and Dinh-Sibony \cite[Example 4.5]{DS},
for every $n \ge 2$, there is an $n$-dimensional
abelian variety $A$
with an automorphism group $G$ satisfying $N(G) = \Aut_0(A) \cong A$, and $G/N(G) = \BZZ^{n-1}$.

\item[(3)]
As seen from the proof of Theorem \ref{ThA},
$X \to X'$ is not an isomorphism
only when: $A$ is $E_{\zeta_3}^3$ or the Jacobian of the
Klein plane quartic, and the index-one Galois cover of $X$ is
a Calabi-Yau threefold and a crepant resolution
of the index-one Galois cover of $X'$;
see Claims \ref{c1} and \ref{c2} and \cite[Theorem (3.4)]{OS} for more details.
\end{itemize}
\end{remark}

With the result above for projective threefolds, we would speculate the following:

\begin{question}\label{Q}
{\rm
Let $G \le \Aut(X)$ be an automorphism group on a minimal
projective variety $X$ of dimension $n \ge 3$
such that the representation $G^* := G | \NS_{\BCC}(X)
\cong \BZZ^{\oplus n-1}$
and every non-trivial element of $G^*$ is of positive entropy.
Is $X$ then birational to the quotient
of an abelian variety by a finite group?
}
\end{question}

We need to assume that $\dim X \ge 3$ in Question \ref{Q} because
a general Wehler $K3$ surface $X$ has Picard number two,
so it is not birational to a {\it generalized}
Kummer $K3$ surface (in the sense of T. Katsura) of necessarily Picard number $\ge 17$.
Here a Wehler $K3$ is the complete intersection $X \subset \BPP^2 \times \BPP^2$ of two hypersurfaces
of bidegrees $(1, 1)$ and $(2, 2)$, where the two projections give rise to
involutions $\iota_i$ so that $g := \iota_1 \circ \iota_2 \in \Aut(X)$
is of topological entropy $\log (7 + 4 \sqrt{3}) > 0$; I thank Serge Cantat for reminding
me about the Wehler surfaces.

Below is an
answer to Question \ref{Q}
under the condition (1) or (2). The minimality of
$(X, G)$ in (1) seems achievable, and necessary too
because the blowup of a
variety with only canonical singularities
contains a rational curve and hence has no \'etale torus covering.

\begin{theorem}\label{ThB}
Let $G \le \Aut(X)$ be an automorphism group on a smooth minimal
projective variety $X$ of dimension $n \ge 3$
such that the representation $G^* := G | \NS_{\BCC}(X)
\cong \BZZ^{\oplus r}$
for some $r = r(G) \ge n-1$,
and every non-trivial element of $G^*$ is of positive topological entropy.
Suppose {\rm either} one of the following two conditions.
\begin{itemize}
\item[(1)]
The pair $(X, G)$ is minimal in the sense that every $G$-equivariant
birational morphism $X \to X'$ onto some $X'$ with only isolated
canonical singularities $($and no other singularities$)$, is an isomorphism.
\item[(2)]
$X$ has no $G$-periodic uniruled subvariety $S$ $(\ne \pt, X)$.
Here, a subvariety is $G$-periodic if it is stabilized by a finite-index
subgroup of $G$.
\end{itemize}
Then $r = n-1$, and there is a $G$-equivariant finite \'etale Galois cover $\tau: A \to X$
for an $n$-dimensional abelian variety $A$, with $\deg \tau \ge 2$ occurring only when
$q(X) = 0$.
\end{theorem}

For the $G$ in Theorem \ref{ThA}(2) or \ref{ThB},
we obtain the following effective
lower bound for the first dynamical degree $d_1(g)$ of $g \in G$.
Set $1 < \sqrt{\delta_n} := \min\{|\gamma_f| > 1 \, ; \, \gamma_f$ is an eigenvalue, with the maximum
modulus, of an integral polynomial $f(x)$ of degree $\le 2n \}$.

\begin{corollary}\label{cThB}
With the notation and assumption
in Theorem $\ref{ThA}$ or $\ref{ThB}$, assume that $r(G) = \dim X - 1 = n - 1 \ge 2$.
Then the first dynamical degree of every $g \in G$ of positive topological entropy satisfies:
$\hskip 1.5pc d_1(g) \ge \delta_n$.
\end{corollary}

When $n = 2$, the result parallel to Corollary \ref{cThB} is that
$d_1(g) \ge \sqrt{\delta_{11}}$ since $\rank_{\BZZ}$ $\, H^2(X, \BZZ) \le 22$
for every minimal surface $X$ with $K_X \equiv 0$ (numerically).
However, when $\dim X \ge 3$, no upper bound of $h^2(X, \BCC)$ is known
even for Calabi-Yau threefolds; this is related to the so called Miles Reid's fantasy:
the moduli space of 3-folds $X$ with $K_X \equiv 0$ may nevertheless be irreducible
(like the case of K3 surfaces).
The lower bound in Corollary \ref{cThB} is independent of Reid's fantasy
which is very hard,
and hence meaningful.

When $\dim X = 2$, McMullen
has proved that
$d_1(g) \ge \lambda_{\Lehmer} = 1.17628081\cdots$ (the Lehmer number)
for all $g \in \Aut(X)$ of positive topological entropy
(cf.~e.g.~\cite[p2]{Mc}).

See Dinh-Sibony \cite{DS}, McMullen \cite{Mc}
and Oguiso \cite{Og09} for related results and references.

\par \vskip 0.5pc
We end the introduction with examples $(X, G)$ where $r(G) = \dim X - 1$
(in the notation of
\ref{ThB})
and $X$ is either rationally connected in the sense of
Campana and Koll\'ar-Miyaoka-Mori, or {\it Calabi-Yau} (i.e., $q(X) = 0$, the Kodaira dimension
$\kappa(X) = 0$, and $K_X \sim_{\BQQ} 0$).

I would like to thank Keiji Oguiso for the
discussion about the example below.

\begin{example}
{\rm
Let $E := \BCC/(\BZZ + Z\sqrt{-1})$ and $A := E^n$ ($n \ge 1$).
Then $\mu_4 = \langle \sqrt{-1} \rangle$ acts diagonally
on the $n$-dimensional abelian variety
$A$ with a few isolated fixed points.
Set $X := A/\mu_4$.
As in \cite[Example 4.5]{DS}, a subgroup $G \cong \BZZ^{\oplus n-1}$ of $\SL_n(\BZZ)$
acts on $A$ (and hence on $X$ since the action of $\mu_4$ is diagonal)
such that every non-trivial element of $G$ is of positive entropy on $A$
(and hence on $X$; cf.~\cite[Lemma 2.6]{Z2} or \cite[Lemma A.8]{NZ}).
Thus we have $r(G) = \dim X - 1$.
If $n \ge 4$ (resp. $n \le 3$) then $X$ is a Calabi-Yau variety
(resp. rationally connected);
see also~\cite[\S 2]{Be}, and \cite[Theorem 1.2, Remark 1.3]{Z4}.
}
\end{example}

\section{Proof of Theorems and Corollaries}

\begin{setup}
We use the conventions in Hartshorne's book, and \cite{KM}.

For an abelian variety $A$, we have $\Aut_{\variety}(A) = T \rtimes \Aut_{\group}(A)$,
where $\Aut_{\variety}(A)$ (resp. $\Aut_{\group}(A)$) is the group
of automorphisms of $A$ as a variety (resp. as a group), and $T = T(A) \cong A$ is the group
of translations of $A$.
By $G|X$, we mean $G \le \Aut(X)$.
\end{setup}

\begin{setup}\label{note1}
Theorem \ref{ThA} (1) follows from the argument in \cite[Theorem 1.2, Remark 1.3]{Z4}
with $X$ replaced by its $G$-equivariant resolution,
$H^2(X, \BZZ)$ by $\NS(X)/(\torsion)$ and the K\"ahler cone by the ample cone;
see also \cite[Proposition 4.8, its remark]{KOZ}.
\end{setup}

\begin{setup}\label{step1}
{\bf Proof of Theorem \ref{ThA} (2)}
\end{setup}

Since $X$ is a projective threefold with $K_X$ nef, the Kodaira dimension
$\kappa(X) \ge 0$. By \cite[Lemma 2.11]{Z4},
$\kappa(X) = 0$; see also \ref{note1}.
By the abundance theorem of Kawamata and Miyaoka,
$IK_X \sim 0$ for some minimal $I = I(X) > 0$, called the {\it index} of $X$.
Let
$$\pi : Y = \oplus_{i=0}^{I-1} \,\, \SO_X(-iK_X) \, \longrightarrow \, X$$
be the {\it index-one covering} so that $K_{Y} \sim 0$.
Lift $G$ to $G \cong G|Y$.
Each summand $\SO_X(-iK_X)$ is an eigenspace of $\Gal(Y/X) = \langle \sigma \rangle$
and we may assume that $\sigma$ acts as the multiple by $\zeta_I^i$ on this summand,
where $\zeta_I = \exp(2 \pi \sqrt{-1} /I)$.
So $G = G|Y$ normalize $\langle \sigma \rangle$ (for later use).
The identification $G|X$ with $G|Y$ also identifies $N(G|X)$ with $N(G|Y)$;
see \cite[Proposition 4.8, its remark]{KOZ},
so $G/N(G) = \BZZ^{\oplus 2}$ holds on $X$ and also on $Y$.

\begin{definition}\label{c2def}
As in \cite{SW}, for a normal projective variety $X$
of dimension $n \ge 3$ and with only canonical singularities,
take a resolution $\iota : Z \to X$ crepant in codimension two
(or minimal in codimension two),
e.g., a composite of a terminalization (cf.~\cite[Corollary 1.4.3]{BCHM})
and a desingularization.
Then we define the (multi) {\it linear form} $c_2(X)$ on $N^1(X) \times \cdots \times N^1(X)$
(with $N^1(X) := \NS(X) \otimes_{\BZZ} \BRR$)
as
$$H_1 \dots H_{n-1} . c_2(X) : = \iota^*H_1\cdots \iota^*H_{n-2} . c_2(Z).$$
By Hironaka's equivariant resolution,
when $G \le \Aut(X)$ is given, we may also choose $\iota$ to be $G$-equivariant,
so that $g^*c_2(X) = c_2(X)$ for all $g \in G$.

We remark that in the case where $K_X$ is nef, $c_2(X) = 0$ (as a linear form) holds if and only if
$H^{n-2} . c_2(X) = 0$ for one ample divisor $H$, by Miyaoka's pseudo-effectivity of
$c_2$ of every terminalization of $X$ and since $N^1(X)$ is spanned by ample divisors $H_i$
(and noting that $H - \varepsilon H_i$ is an ample $\BQQ$-divisor for small $\varepsilon$).
\end{definition}

\begin{claim}\label{c1}
\begin{itemize}
\item[(1)]
If $Z \in \{X, Y\}$, 
then $H_Z . c_2(Z) = 0$ for a nef and big $\BRR$-divisor $H_Z$.
\item[(2)]
If $q(X) > 0$ $($resp. $q(Y) > 0)$, then $X$ $($resp. $Y)$ is an abelian variety,
so Theorem $\ref{ThA} (2)$ is true.
\item[(3)]
If $c_2(X) = 0$ or $c_2(Y) = 0$ as linear form, then Theorem $\ref{ThA} (2)$ is true.
\end{itemize}
\end{claim}

We prove Claim \ref{c1} (2) first.
We may assume that the $s$-th group $G^{(s)}$ in the derived series of $G$
satisfies $G^{(s)} | \NS_{\BRR}(X)$ $= \id$.
Take a finite-index subgroup $G_1 \le G$ such that
$G_1 | \NS_{\BCC}(Y)$ is $Z$-connected and $G_1/N(G_1) \cong \BZZ^{\oplus 2}$ still holds
on $X$ and $Y$.
Take an ample Cartier divisor $M \subset X$ and write
$M' = \pi^*M$. Then the class $[M'] \in \NS_{\BRR}(Y)$
is $G_1^{(s)}$-invariant. Now a result of David Lieberman
(and Hironaka's equivariant resolution)
imply that $\Aut_{M'}(Y) := \{g \in \Aut(Y) \, | \, g^*M' \equiv M'$
in $\NS_{\BRR}(Y)\}$ is a finite-index over-group of $\Aut_0(Y)$
(where the latter acts trivially on $\NS_{\BRR}(Y)$); see \cite[Lemma 2.23]{Z2}.
So $G_1^{(s)} | \NS_{\BCC}(Y)$ is a finite group and hence is trivial
because $G_1^{(s)} | \NS_{\BCC}(Y)$ is also $Z$-connected.
Thus, $G_1 | \NS_{\BCC}(Y)$ is solvable and also $Z$-connected
with $r(G_1|Y) = 2$.

By \cite[Lemma 2.13]{Z4} (see also \ref{note1}),
either $q(Y) = 0$, or the albanese map $\alb_Y : Y \to A := \Alb(Y)$
is a well defined morphism (for $Y$ having only rational singularities and using
\cite[Lemma 8.1]{Ka})
and is birational and surjective.
In the latter case, $Y \cong A$ (and Theorem \ref{ThA} (2) is true as
in (3) below),
because $K_{Y} \equiv 0 \equiv K_A$, and $\alb_Y$
is neither a small contraction nor crepant
since $A$ is smooth.
The same argument applies to $X$. This proves (2).

(1) We treat only $Y$, since the case $X$ is similar and simpler.
As in (2), replacing $G$ with its finite-index subgroup,
we may assume that $G | \NS_{\BCC}(Y)$ is solvable and Z-connected.

We use the argument in the proof of \cite[Lemma 5.2]{KOZ}.
In particular, there is a common nef eigenvector
$0 \ne L \in \NS_{\BRR}(Y) \cap c_2(Y)^{\perp}$ for $G$ on $Y$.
For $g \in G$, write $g^*L \equiv \chi(g) L$ with $\chi(g) \in \BRR_{>0}$.
Consider the homomorphism below
$$f : G | Y \longrightarrow \BRR, \hskip 2pc g \mapsto \log \chi(g).$$
We have $N(G) \subseteq \Ker(f)$.
Consider the case
$N(G) \ne \Ker(f)$. Take $g \in \Ker(f) \setminus N(G)$.
By the generalized Perron-Frobenius theorem \cite{Bi},
there are nef $\BRR$-divisors $L_g^{\pm}$ such that
$(g^{\pm})^* L_{g}^{\pm} \equiv d_1(g^{\pm}) L_g^{\pm}$.
Here $d_1(h)$ denotes the first dynamical degree of $h \in \Aut(Y)$.
We have $L . L_g^+ . L_g^- \ne 0$ by \cite[Lemma 4.4]{DS},
so $H_Y := L + L_g^+ + L_g^-$ is nef and big with $H_Y . c_2(Y) = 0$.
Indeed, $H_Y^3 \ge L . L_g^+ . L_g^- > 0$; also
$L . c_2(Y) = 0$ by the choice of $L$, and
$L_g^{\pm} . c_2(Y) = (g^{\pm})^*L_g^{\pm} . (g^{\pm})^*c_2(Y) =
d_1(g^{\pm}) L_g^{\pm} . c_2(Y)$, so $L_g^{\pm} . c_2(Y) = 0$, since
$d_1(g^{\pm}) > 1$ for $g \not\in N(G)$.

Thus we may assume that $\Ker(f) = N(G)$.
If every $\chi(g)$ is $1$, or $d_1(g)$, or $1/d_1(g^{-1})$, then
$\Imm(f)$ is discrete in $\BRR$ by \cite[Corollary 2.2]{DS};
so $\BZZ^{\oplus 2} \cong G/N(G) \cong \Imm(f) = \BZZ^{\oplus s}$ with $s \le 1$,
absurd. Therefore, we may assume that
$\chi(g) \ne 1$, $d_1(g)$, or $1/d_1(g^{-1})$ for
some $g \in G$. Then we have $L . L_g^+ . L_g^- \ne 0$ by \cite[Lemma 4.4]{DS}.
Thus Claim \ref{c1} (1) is true as above by taking $H_Y := L + L_g^+ + L_g^-$.

(3) If $c_2(Z) = 0$ as a linear form for some $Z \in \{X, Y\}$,
then $Z = A/H$ for an abelian variety $A$ and $H \le \Aut_{\variety}(A)$
acting on $A$ freely in codimension two (cf.~\cite[Cor, p.266]{SW}).
Replacing $A$, we may assume 
$A \to X$ is (Galois and)
the Albanese
closure in codimension one in the sense of \cite[Lemma 2.12]{nz2}
(Kawamata's characterization
of abelian variety may also be used), so that $G$ lifts to $G \cong G|A \le \Aut_{\variety}(A)$.
This proves (3) and
Claim \ref{c1}.

\par \vskip 1pc
We resume the proof of Theorem \ref{ThA} (2).
By Claim \ref{c1}, we may assume that $q(Z) = 0$ and $c_2(Z) \ne 0$
for $Z = X$, $Y$.
By \cite[Theorem (3.4)]{OS}, $Y$ is smooth and there are a
birational crepant morphism $Y \to Y'$ (the unique contraction of all divisors
perpendicular to $c_2(Y)$) and an \'etale-in-codimension-two Galois covering
$A \to Y'$ from an abelian variety $A$ of dimension three.
Further, either $A = E_{\zeta_3}^3$ $($the product of three copies of the elliptic curve
of period $\zeta_3 = \exp(2 \pi \sqrt{-1} /3))$ and $|\Gal(A/Y')|$ divides $27$,
or $A$ is the Jacobian of the Klein plane quartic $\{X_0X_1^3 + X_1X_2^3 + X_2X_0^3 = 0\}$
and $\Gal(A/Y') \cong \BZZ/(7)$.
Note that $(G|Y) . \langle \sigma \rangle$ (with $\langle \sigma \rangle = \Gal(Y/X)$
normalized by $G|Y$ as mentioned early on)
descends to a regular action
on $Y'$ and we let $Y' \to X' := Y'/\langle \sigma \rangle$ be
the quotient map. Then there is a birational morphism $X \to X'$
such that the two natural composites below coincide:
$Y \to Y' \to X'$,
$Y \to X \to X'$.
Our $G | X$ (and $G|Y'$) descend to a regular action $G | X' \cong G$.
Since $Y \to X$ is \'etale in codimension one,
so is $Y' \to X'$;
further, $X \to X'$ is crepant because
so is $Y \to Y'$.
Therefore, $X'$ has only canonical
singularities (and $\BQQ$-factorial because $A \to X'$ is finite)
with $K_{X'} \sim_{\BQQ} 0$.
Note that the composite $A \to Y' \to X'$ is \'etale in codimension one.
Replacing $A$, we may assume that $A \to X'$ is (Galois and) the Albanese
closure in codimension one in the sense of \cite[Lemma 2.12]{nz2},
so that $G \cong G|X'$ lifts to $G \cong G|A \le \Aut_{\variety}(A)$.
This proves Theorem \ref{ThA} (2).

\begin{setup}\label{step2}
{\bf Proof of Theorem \ref{ThA} (3)}
\end{setup}

The identification of $G|X$, $G|X'$ and $G|A$ also identifies $N(G|X)$, $N(G|X')$
and $N(G|A)$;
see \cite[Proposition 4.8, its remark]{KOZ},
so $G/N(G) = \BZZ^{\oplus 2}$ holds on $X$, $X'$ and $A$.
Let $G_0$ be a finite-index normal subgroup of $G$ such that $G_0 | \NS_{\BCC}(A)$ is
$Z$-connected.
Since $\NS_{\BCC}(X')$ pulls back to a subspace of $\NS_{\BCC}(A)$, our
$G_0 | \NS_{\BCC}(X')$ is also $Z$-connected.
Further, $G_0 | \NS_{\BCC}(X)$ is $Z$-connected,
since $\NS_{\BCC}(X)$ is spanned by the pullback of $\NS_{\BCC}(X')$
and finitely many irreducible components in the exceptional locus of $X \to X'$
which is stable under the actions of $G|X$ and $G_0|X$.
Note that $N_0 := G_0 \cap N(G)$ equals $N(G_0)$ on all of $X$, $X'$ and $A$.
We also have $G_0/N_0 \cong \BZZ^{\oplus 2}$ and $r(G_0) = 2$.

\begin{claim} \label{c3}
$G_0 | \NS_{\BCC}(A)$ is solvable.
If $N_0$ is finite, then $N(G)$ is finite.
\end{claim}

We prove Claim \ref{c3}.
The first part follows from the proof of Claim \ref{c1} and the solvability
of the action $G_0$ ($\le G$) on $\NS_{\BCC}(X)$ and hence on $\NS_{\BCC}(X')$.
Suppose $s_2 := |N_0| < \infty$. Set $s_1 := |G:G_0|$ and $s = s_1s_2$.
Take any $n \in N(G)$. Then $n^{s_1} \in G_0 \cap N(G) = N_0$ and hence
$n^s = \id$. Thus $N^* := N(G) | \NS_{\BCC}(X)$ is a periodic group with bounded exponent,
so it is a finite group by Burnside's theorem.
Hence $N(G) \le \Aut_M(X)$, where $M := \sum_{n^* \in N^*} n^*M'$ with an ample divisor $M' \subset X$.
Now $N(G)$ is finite by the assumption $\Aut_0(X) = (1)$ and Lieberman's result as in Claim \ref{c1}.
This proves Claim \ref{c3}.

\vskip 1pc
We resume the proof of Theorem \ref{ThA} (3).
By Claim \ref{c3} and replacing $G$ with $G_0$, we may assume that $G|\NS_{\BCC}(A)$ is already $Z$-connected
(and also solvable).
Let $\bar{G}$ and $\bar{N}$ be the images of $G$ and $N(G)$ under the projection
$\Aut_{\variety}(A) = T \rtimes \Aut_{\group}(A) \to \Aut_{\group}(A)$,
where $T = \Aut_0(A)$ is the group of translations.

We are going to apply Oguiso \cite[Lemma 2.5]{Og}.
Consider the faithful matrix representation $\bar{G} | H^0(A, \Omega_A^1)$.
We say that $g \in \bar{G}$ is {\it unipotent} if so is its representation matrix,
and let $U$ be the set of unipotent elements in $\bar{G}$.
An element $g$ of $\bar{G}$ is in $\bar{N}$ if and only if all eigenvalues
of its matrix representation are of modulus $1$, i.e., they are of roots of $1$
by Kronecker's theorem (and with bounded minimal polynomial over $\BQQ$).
Thus $U \subseteq \bar{N}$ and there is an $s > 0$ such that $n^s \in U$ for all $n \in \bar{N}$.
If every element of $\bar{N}$ is periodic then these periods divide $s$;
thus, by Burnside's theorem,
$\bar{N} | H^{1,1}(A)$
is a finite group, so $N(G) | H^{1,1}(A)$ and hence
$N(G) | \NS_{\BCC}(X)$
(embedded in the former by the pullback)
and even $N(G)$ are all finite groups (cf.~the proof of Claim \ref{c3}).

Therefore, we may assume that $U$ contains an element of infinite order.
In particular, the pointwise fixed set $B := A^U = \{a \in A \, | \, u(a) = a, \,\, \forall u \in U\}$
satisfies $0 \le B < A$.
By \cite[Lemma 2.5]{Og}, $U$ is a normal subgroup of $\bar{N}$ (and also of $\bar{G}$
because conjugate action preserves the unipotent property).
It is known then that $U|H^0(A, \Omega_A^1)$ can be regarded as a subgroup of the unipotent
group $U(q(A), \BCC)$ and of the upper triangular group $T(q(A), \BCC)$,
so the matrix representation of $U$ has a common eigenvector
(corresponding to the unique eigenvalue $1$).
Thus $B = \cap_{u \in U} \, \Ker(u - \id)$ satisfies $0 < B < A$.
Since $U$ is normal in $\bar{G}$, our subtorus $B$ is $\bar{G}$-stable.
So $G$ permutes the cosets $A/B$. Thus we have an induced
$G$-equivariant fibration $A \to A/B$. Therefore, $2 = r(G) \le \dim A - 2 = 1$
by the proof of \cite[Lemma 2.10]{Z4} (see also \ref{note1}).
This is a contradiction.
The proves Theorem \ref{ThA} (3).

\begin{setup}
{\bf Proof of Corollary \ref{CorB}}
\end{setup}

Let $G_1 \le G$ be the inverse of the identity connected component for the Zariski-closure
of $G^*$ in $\GL(\NS_{\BCC}(X))$. Then $G_1$ satisfies the hypothesis of Theorem \ref{ThA}
and $r(G_1) = r(G) = 2$, since $|G:G_1| < \infty$ and $\Ker (G \to G^*) = N(G)$
by the assumption.
Thus Corollary \ref{CorB} holds. Indeed,
the $G$-equivariant property is true because the morphisms involved in the proof of
Theorem \ref{ThA},
like the index-one cover, or Albanese-closure in codimension one, or the $c_2$-birational contraction
in \cite[Theorem (3.4)]{OS} are all canonical.

\begin{setup}
{\bf Proof of Theorem \ref{ThB}}
\end{setup}

At first, we assume no smoothness of $X$, but assume
$X$ has only canonical singularities.
Note that $r = r(G) = n-1$ (cf.~\ref{note1}, \cite[Theorem 4.7]{DS}, \cite[Theorem 1.2, Remark 1.3]{Z4}).

Applying \cite[Theorems 4.7 and 4.3]{DS} to a $G$-equivariant resolution of $X$,
there are nef $\BRR$-divisors $L_1, \dots, L_n$ as common eigenvectors
of $G$ such that $L_1 \cdots L_n \ne 0$ and the homomorphism
below is an isomorphism onto a spanning lattice
(where we write $g^*L_i \equiv \chi_i(g) L_i$):
$$f : \, G | \NS_{\BRR}(X) \longrightarrow \BRR^{n-1}, \hskip 2pc g \mapsto (\log \chi_1(g), \dots, \log \chi_{n-1}(g)).$$
Since $L_1 \dots L_n = g^*(L_1 \cdots L_n) = \chi_1(g) \cdots \chi_n(g) L_1 \dots L_n$
we have $\chi_1 \dots \chi_n = 1$.
Set $H := \sum_{i=1}^n L_i$. Since $H^n \ge L_1 \dots L_n > 0$, our $H$ is nef and big.

\begin{claim}\label{c4}
$H^{n-1} . K_X = 0$, and $H^{n-2} . c_2(X) = 0$.
\end{claim}

We prove Claim \ref{c4}.
Take $M := L_1^{i_1} \cdots L_{n}^{i_n}$ with $\sum_{k=1}^n i_k = s$.
When calculating $M . c_p(X)$ for $p = 1$ or $p = 2$,
we let $s = n-p$.
For $g \in G$, we have $g^*M = e(g) M$ with $e(g) = \chi_1(g)^{i_1} \cdots \chi_n(g)^{i_n}$.
Since $M . c_p(X) = g^*M . g^*c_p(X) = e(g) M . c_p(X)$ and since $H^{s}$ is a combination of such $M$,
it suffices to show that $e(g) \ne 1$ for some $g \in G$ (so that $M . c_p(X) = 0$).
Suppose the contrary that $e(g) = 1$ for all $g \in G$.
Taking $\log$ and using $\chi_1 \dots \chi_n = 1$, we have
$(i_1-i_n) \log \chi_1 + \cdots + (i_{n-1}-i_n) \log \chi_{n-1} = 0$ on $G$.
Since the image of the homomorphism $f$ above is a spanning lattice,
this happens only when $i_1 - i_n = \cdots = i_{n-1} - i_n = 0$.
Thus $n-1 \ge s = \sum_{k=1}^n i_k = ni_1$, so $i_1 = 0$ and hence $s = 0$, absurd.
This proves Claim \ref{c4}.

\vskip 1pc
Since $H^{n-1} . K_X = 0$ and $K_X$ is nef, we have $K_X \equiv 0$ by
\cite[Lemma 2.2]{nz2}. So $K_X \sim_{\BQQ} 0$ by \cite[Theorem 8.2]{Ka}.
By the $\BRR$-divisor version \cite[Theorem 3.9.1]{BCHM}
of Kawamata's base point freeness theorem,
there is a birational morphism $\gamma : X \to X'$ such that
$H = \gamma^* H'$ for an ample $\BRR$-divisor $H'$.
The result below should be well known, but we prove it
for the lack of reference. As said early,
every singularity of $X$ is assumed to be canonical only.

\begin{claim}\label{c2}
$\gamma : X \to X'$ is crepant, and hence
$X'$ has only canonical singularities.
Further, the indices of $X$ and $X'$ coincide and are denoted as $J$
so that $JK_X \sim 0$ and $JK_{X'} \sim 0$.
\end{claim}

We prove Claim \ref{c2}.
Let $I := I(X)$ be the index of $X$ and
$$\pi : \, Y = \Spec \oplus_{i=0}^{I-1} \,\, \SO_X(-iK_X) \, \longrightarrow \, X$$
the index-one cover which is \'etale in codimension one, where $K_Y \sim 0$.
Let $Y \to Y' \to X'$ be the Stein factorization
of $Y \to X \to X'$. Then $Y' \to X'$ is \'etale in codimension one
because so is $Y \to X$. By \cite[Proposition 5.20]{KM},
every singularity of $Y$ is canonical and hence rational.
Applying Koll\'ar's torsion freeness result for higher direct image of dualizing
sheaf to the {\it composite} of $Y \to Y'$
and a resolution of $Y$, the first condition in \cite[Proposition 3.12]{Ko}
is satisfied by this composite because $\SO(K_Y) \cong \SO_Y$ now
(cf.~also \cite[Corollary 2.68]{KM}),
so $Y'$ has only rational singularities. Applying \cite[Lemma 5.12]{KM}
to the above composite, we have $\SO_{Y'}(K_{Y'}) \cong \SO_{Y'}$.
Thus $Y'$ has only Gorenstein canonical singularities,
so $X'$ has only log terminal singularities by \cite[Proposition 5.20]{KM},
and $IK_{X'} \sim 0$;
further, $Y \to Y'$ is crepant because $K_Y \sim 0 \sim K_{Y'}$.
The commutative diagram in the Stein factorization above then implies
that $X \to X'$ is crepant. The last part of Claim \ref{c2}
is true because $I(X')K_{X'} \sim 0$ implies $I(X')K_X \sim 0$
by considering the fibre product over $X'$, of $X \to X'$ and the index-one covering of $X'$.
This proves Claim \ref{c2}.

\vskip 1pc
$(H')^{n-2} . c_2(X') = H^{n-2} . c_2(X) = 0$ implies $c_2(X') = 0$
as a linear form (cf.~\ref{c2def}). Our
$\gamma:$ $X \to X'$ is $G$-equivariant, because a curve $C \subset X$
is contracted by $\gamma$ if and only if $H . C = 0$; if and only if
$L_i . C = 0$ (for all $i$, since $L_i$ is nef); if and only if $L_i . g_*C = 0$ (since $L_i$ is
semi $G$-invariant); if and only if $g_*C$ is contracted by $\gamma$.
This argument and \cite[Remark (2), p46]{KM} show that $L_i = \gamma^*L_i'$ for
some nef $L_i'$ on $X'$ (with $g^*L_i' \equiv \chi_i(g) L_i'$). Thus, $H' = \sum_{i=1}^n L_i'$.

We assert that $X'$ has no $G$-periodic subvariety $S'$ of
dimension $s \in \{1, \dots, n-1\}$. Indeed, the proof of Claim \ref{c4} applied to a finite-index
subgroup of $G | X'$, would show that $(H')^s . S' = 0$, contradicting the ampleness of $H'$.

If $\gamma$ is not an isomorphism, then the exceptional locus $\Exc(\gamma)$ is non-empty
whose irreducible components $E_i$ are permutated by $G$ and hence
$G$-periodic. So $\gamma(E_i)$ are $G$-periodic and hence a point by the assertion above.
Thus $E_i$ is covered by fibres
of $\gamma$ and hence uniruled, since every fibre of a partial resolution
of a log terminal singularity is rationally chain connected by Hacon-McKernan's
solution to Shokurov's conjecture.

We remark that for every $G$-periodic subvariety $\pt \ne S \subset X$,
the image $\gamma(S) \subset X'$ is $G$-periodic
and hence a point by the assertion above, thus $S \subseteq \Exc(\gamma)$.

By either conditoin,
we may assume
$\gamma = \id$.
As in
Theorem \ref{ThA}(2), Theorem \ref{ThB} follows from
the vanishing of $c_i(X)$ ($i = 1, 2$),
Bogomolov decomposition \cite{Be} and {\it smoothness} of $X$.

\begin{setup}
{\bf Proof of Corollary \ref{cThB}}

It follows from that $\rank_{\BZZ} \, H^1(A, \BZZ) = 2n$ and $d_1(g) = d_1(g_A)$
(cf.~\cite[Lemma 2.6]{Z2} or \cite[Lemma A.8]{NZ});
here $g_A$ on $A$ is the lifting of $g$.

\end{setup}

\end{document}